\documentclass[final,leqno,onefignum,onetabnum]{siamltex1213}
\usepackage{amsmath,amssymb}
\usepackage{rotating}
\usepackage{booktabs}
\usepackage{multirow}

\title{On the complexity of local search in unconstrained quadratic binary optimization} 

\author{D{\'a}vid Papp\thanks{Department of Mathematics, North Carolina State University, Raleigh, NC, USA. Email: dpapp@ncsu.edu}}

\date{\today}

\newtheorem{thm}{Theorem}

\newcommand{\deletethis}[1]{{}}

\newcommand{\vb}{\mathbf{b}}

\newcommand{\vQ}{\mathbf{Q}}

\newcommand{\vx}{\mathbf{x}}

\newcommand{\T}{\mathrm{T}}

\newcommand{\defeq}{\ensuremath{\overset{\mathrm{def}}{=}}}

\begin{document}
	\maketitle
	\slugger{siopt}{xxxx}{xx}{x}{x--x}
	
\begin{abstract}
We consider the problem of finding a local minimum of a binary quadratic function, and show by an elementary construction that every descending local search algorithm takes exponential time in the worst case.
\end{abstract}
	
	
	\begin{AMS}\end{AMS}
	
	\pagestyle{myheadings}
	\thispagestyle{plain}
	\markboth{D{\'A}VID PAPP}{LOCAL SEARCH FOR BINARY QUADRATIC OPTIMIZATION IS EXPONENTIAL}
	
\section{Introduction}
We consider the problem of finding a local minimum of a multivariate quadratic polynomial with binary variables, that is, a local optimal solution to
\begin{equation}\label{eq:1} \min_{\vx\in\{0,1\}^n} f(\vx) \defeq \tfrac{1}{2}\vx^\T\vQ\vx + \vb^\T\vx,
\end{equation}
where $\vQ$ is a given $n\times n$ integer symmetric matrix and $\vb$ is a given integer $n$-vector. Local minima is defined with respect to the usual ``$k$-flip'' neighborhood for a fixed integer $k$ that is independent of the problem instance; that is, a vector $\vx$ is a \emph{local minimum of $f$} if $f(\vx)$ cannot be decreased by changing up to $k$ components of $\vx$ (from 0 to 1, or vice versa). A \emph{local search method} starts from a given initial binary vector (that is part of the input, along with $\vQ$ and $\vb$ that define $f$), and in every iteration decides whether the current vector is a local minimum or not.  If the answer is positive, it returns the current binary vector; if the answer is negative, it moves to a better binary vector by flipping up to  $k$ components.

Finding the global optimal solution to \eqref{eq:1} is widely known to be NP-hard. Similarly, binary quadratic optimization, and even its special case, the maximum cut problem is PLS-complete \cite{Yannakakis2009}. The complexity class PLS contains every optimization problem satisfying the properties that (i) a feasible solution can be computed in polynomial time, (ii) the objective function is computable in polynomial time, and (iii) there is a polynomial time algorithm that, for every feasible solution, decides whether the solution is locally optimal, and if the answer is negative, computes a better neighbor. In this definition, \emph{polynomial time} refers to polynomial in the size of the input instance, which is in turn equivalent to polynomial in $n$ and in $\log U$, where $U$ is the absolute value of the largest component in $\vQ$ and $\vb$.

Just as the NP-hardness of global optimization in \eqref{eq:1} implies that a polynomial-time algorithm for the solution of the global optimization problem would yield a polynomial time solution to every problem in NP, the implication of the PLS-completeness of local optimization is that a polynomial time algorithm for finding a local optimal solution to \eqref{eq:1} would yield a polynomial time solution to every problem in PLS. It is also known that there are several combinatorial optimization problems for which every local search method (with the appropriate definition of the neighborhood) is exponential. For example, \cite{JohnsonPapadimitriouYannakakis1988} and \cite{Krentel1989} have proved that for every sufficiently large $k$, there are instances of the Traveling Salesman Problem (TSP) and initial tours such that, starting with the initial tour, every sequence of improving $k$-changes that terminates in a $k$-optimal tour has exponential length; \cite{ChandraKarloffTovey1999} has further improved this by showing the same for every $k\geq 3$. (An improving $k$-change in TSP consists of removing no more than $k$ edges from a tour and adding as many distinct edges as the number of removed edges to form a new tour.)

While none of the above results have any implications on the complexity of finding a local optimum for \eqref{eq:1} using algorithms that are not local search methods, they indicate that there are ``hard'' instances of this problem for which \emph{every} local search method takes exponential time in the worst case. 

The goal of this short note is to provide an explicit, compact construction of such an instance. This also provides a counterexample to \cite[Theorems 2.15 and 2.21]{Chen2015}, which state, incorrectly, that a particular $1$-flip (respectively, $2$-flip) local search method terminates in no more than $n$ iterations (resp., $O(n^2)$ iterations) for every $n$-variable problem. In \mbox{Section 2}, only the 1-flip neighborhood is considered. In \mbox{Section 3}, it is shown that hard instances for the 1-flip neighborhood can be extended to similar instances for the $k$-flip neighborhood, for any fixed value of $k$.

\section{The construction for the 1-flip neighborhood}

Define, for $n=1,5,9,\dots$, the sequence of $n$-variate polynomials $f_n$ recursively as follows, using the shorthand $\vx$ for the vector $(x_1,\dots,x_n)$:
\begin{multline*}
f_1(x_1) = -x_1,\\
f_{n+4}(\vx,x_{n+1},x_{n+2},x_{n+3},x_{n+4}) = f_n(\vx) - M_n \Big( (-2n+1) x_{n+1} - 2 x_{n+2} -2 x_{n+3} - 6n x_{n+4}\\
 + \left(\sum_{i=1}^n x_i\right) (2 x_{n+1} - 4 x_{n+2} - 2 x_{n+3} + 4 x_{n+4})\\
 + 7n x_{n+1}x_{n+2} + 3 x_{n+2}x_{n+3} + 7n x_{n+3}x_{n+4} \Big),
\end{multline*}
with
\begin{equation}\label{eq:Mn}
M_n = \max_{\vx\in\{0,1\}^n} f_n(\vx) - \min_{\vx\in\{0,1\}^n} f_n(\vx) + 1.
\end{equation}
For example,
\begin{multline*}
f_5(x_1,\dots,x_5)= -x_1+2 x_2+4 x_3+4 x_4+12 x_5\\
-4x_1x_2 +8x_1x_3 +4x_1x_4 -8x_1x_5 -14x_2x_3  -6x_3 x_4 -14x_4x_5.
\end{multline*}

We claim that every local search method starting from $\vx=(0,\dots,0)$ takes exponentially many steps to reach the local optimum (with respect to the 1-flip neighborhood) $\vx=(1,\dots,1)$. 

\begin{thm}\label{thm:1}
Every 1-flip local search method applied to $f_n$ starting from $\vx=(0,\dots,0)$ takes exponentially many steps in the size of $f_n$ to reach the local minimum $\vx=(1,\dots,1)$.
\end{thm}
\begin{proof}
Throughout the proof, we shall use the shorthand $S_n = \sum_{i=1}^n x_i$.

First, we show that the size of $f_n$ is polynomial in $n$. This is equivalent to saying that $\log(M_n)$ can be bounded from above by a polynomial in $n$. The maximum and minimum of $f_n(\vx)$ over binary $n$-vectors can be bounded from above and below, respectively, by the sum of the positive (resp., negative) coefficients of $f_n$. Therefore, we have $M_1=2$ and using $S_n\leq n$ and $1\leq n$, we obtain
\[M_{n+4} \leq M_n + M_n( (8n+6S_n+3) + (14n+6S_n+3) ) \leq 41 n M_n. \]
It follows that $M_n \leq 41^n n!$, therefore $\log(M_n)$ is of order $O(n\log n)$.

To prove our theorem, it suffices to show that every local search method starting from $(0,\dots,0)$ terminates at $(1,\dots,1)$, and that it does so after exponentially many iterations in $n$. To keep the case analysis that follows as clear as possible, we use the observation that for quadratic binary functions, the partial derivatives are the same as partial discrete derivatives:
\begin{equation}\label{eq:derivatives}
\frac{\partial f}{\partial x_i}(\vx) = f(x_1,\dots,x_{i-1},1,x_{i+1}\dots,x_n)-f(x_1,\dots,x_{i-1},0,x_{i+1}\dots,x_n).
\end{equation}
The identity \eqref{eq:derivatives} implies that changing a variable from zero to one results in a decrease in the function value if and only if the corresponding partial derivative is negative. 

Clearly, for $n=1$, every local search starting at $0$ will move to $1$ in a single step.

The inductive step is summarized in Table \ref{tbl:1}, which shows the partial derivatives of $f_{n+4}$ and their signs at various points. From the table, we can infer the paths followed by local search methods starting at zero:
\begin{enumerate}
\item As long as $S_n < n$, that is, until $(x_1,\dots,x_n)=(1,\dots,1)$, none of the last four variables can be flipped to $1$. Therefore, by induction, all local search methods applied to $f_{n+4}$ will initially follow a decreasing path from $(0,\dots,0)$ to $\vx^{(1)}=(1,\dots,1,0,0,0,0)$.

\item At $\vx^{(1)}$, the only $0$ component with a corresponding negative partial derivative is the one with respect to $x_{n+1}$, and no $1$ component has a corresponding positive partial derivative. Therefore $x_{n+1}$ must to flipped to $1$, leading to $\vx^{(2)}=(1,\dots,1,1,0,0,0)$.

\item Similarly, at $\vx^{(2)}$, $x_{n+2}$ must be flipped to $1$, leading to $\vx^{(3)}=(1,\dots,1,1,1,0,0)$.

\item At $\vx^{(3)}$, and also at every point of the form $(\vx,1,1,0,0)$ with $S_n>0$, the last four components cannot be flipped, but flipping any of the first $n$ variables leads to a decrease in the objective value. Therefore, in the next $n$ steps, every local search method will flip back to zero each of the variables $x_1$ through $x_n$ in some order, leading to $\vx^{(4)}=(0,\dots,0,1,1,0,0)$.

\item At $\vx^{(4)}$, $x_{n+3}$ must be flipped, leading to $\vx^{(5)}=(0,\dots,0,1,1,1,0)$.

\item At $\vx^{(5)}$, $x_{n+4}$ must be flipped, leading to $\vx^{(6)}=(0,\dots,0,1,1,1,1)$.

\item Once $x_{n+1}=\cdots=x_{n+4}=1$, all of the last four partial derivatives remain negative, and these variables cannot be flipped back to $0$, regardless of what values the variables $x_1$ through $x_n$ take. Therefore, continuing from $\vx^{(6)}$, local search methods will only manipulate the first $n$ variables, following a decreasing path that is a translation of a decreasing local search path for $f_n$ starting at zero, leading, by induction, to the point $(1,\dots,1)$.

\item At $(1,\dots,1)$, all partial derivatives are negative; this is a local minimum.
\end{enumerate}

We conclude that every local search method applied to $f_{n+4}$, starting from $(0,\dots,0)$ will terminate at $(1,\dots,1)$, and they do so after more than twice as many steps as any method applied to $f_n$ takes. Therefore, the number of steps until termination is of order $O(2^{n/4})$.

More precisely, if $s_n$ denotes the number of steps for $f_n$, then we have
\[s_1 = 1 \text{ and } s_{n+4} = 2s_n+n+4 \text{ for }n=1,5,9,\dots,\]
which yields (by induction) $s_n=5\cdot 2^{\frac{n+3}{4}}-n-8$.
\end{proof}

\begin{table}[!t]
\begin{tabular}{lccccc}
\toprule
Partials & at $(\vx,0,0,0,0)$ & at $(\vx,1,0,0,0)$ & at $(\vx,1,1,0,0)$ & at $(\vx,1,1,1,0)$ & at $(\vx,1,1,1,1)$\\
\midrule
\multirow{2}{*}{$\frac{1}{M_n}\frac{\partial f_{n+4}}{\partial x_i}$} & \multirow{2}{*}{$\frac{1}{M_n}\frac{\partial f_{n}}{\partial x_i}$} & $-2+\frac{1}{M_n}\frac{\partial f_{n}}{\partial x_i}$ & $2+\frac{1}{M_n}\frac{\partial f_{n}}{\partial x_i}$ & $4+\frac{1}{M_n}\frac{\partial f_{n}}{\partial x_i}$ & \multirow{2}{*}{$\frac{1}{M_n}\frac{\partial f_{n}}{\partial x_i}$}\\
& & $<0$ & $>0$ & $>0$\\
\midrule
\multirow{2}{*}{$\frac{1}{M_n}\frac{\partial f_{n+4}}{\partial x_{n+1}}$} & $-1+2n-2S_n$ & $-1+2n-2S_n$  & $-1-5n-2S_n$  & $-1-5n-2S_n$ & $-1-5n-2S_n$ \\
 & $>0$ iff $S_n\neq n$ & $>0$ iff $S_n\neq n$ & $<0$ & $<0$ & $<0$\\
\midrule
\multirow{2}{*}{$\frac{1}{M_n}\frac{\partial f_{n+4}}{\partial x_{n+2}}$} & $2+4S_n$   & $2-7n+4S_n$ & $2-7n+4S_n$ & $-1-7n+4S_n$ & $-1-7n+4S_n$\\
 & $>0$ & $<0$ & $<0$ & $<0$ & $<0$ \\
\midrule
\multirow{2}{*}{$\frac{1}{M_n}\frac{\partial f_{n+4}}{\partial x_{n+3}}$} & $2+2S_n$   & $2+2S_n$    & $-1+2S_n$     & $-1+2S_n$    & $-1-7n+2S_n$ \\
 & $>0$ & $>0$ & $>0$ iff $S_n\neq0$ & $>0$ iff $S_n\neq0$ & $<0$\\
\midrule
\multirow{2}{*}{$\frac{1}{M_n}\frac{\partial f_{n+4}}{\partial x_{n+4}}$} & $6n-4S_n$  & $6n-4S_n$   & $6n-4S_n$   & $-n-4S_n$    & $-n-4S_n$ \\
 & $>0$ & $>0$ & $>0$ & $<0$ & $<0$\\
\bottomrule
\end{tabular}
\caption{The partial derivatives of $f_{n+4}$ and their signs at various points. $S_n$ denotes $\sum_{i=1}^n x_i$. A partial derivative is negative if and only if setting the corresponding variable to $1$ (without changing the other variables) yields a lower function value than setting the corresponding variable to $0$.}\label{tbl:1}
\end{table}

\section{Extension to the $k$-flip neighborhood}

Let $k \geq 2$ be a fixed integer, and define the $kn$-variate polynomials $g^{(k)}_n$ as
\[g^{(k)}_n(x_{11}, x_{12}, \dots, x_{1n}, x_{21}, \dots, x_{kn}) = f_n(x_{11}, \dots, x_{1n}) + M_n \sum_{i=2}^k \sum_{j=1}^n (x_{ij}-x_{1j})^2, \]
with the same $M_n$ as defined in \eqref{eq:Mn}. The size of $g^{(k)}_n$ is polynomial in $n$, for every fixed value of $k$.

In essence, $g^{(k)}_n$ has $k$ copies ($x_{1j}$ through $x_{kj}$) of each variable $x_{1j}$ of $f_n$, and attains the same values as $f_n$ whenever these copies are equal. At points where any two copies of the same variable are different, $g^{(k)}_n$ takes a strictly higher value than the maximum of $f$.

Let us call a vector $\vx=(x_{11},\dots,x_{kn})$ \emph{simple} if it satisfies $x_{ij}={x_{1j}}$ for every $i=2,\dots,k$ and $j=1,\dots,n$. Based on the previous paragraph, every simple vector is a local minimum of $g_n^{(k)}$ within the $\ell$-flip neighborhoods for all $\ell=1,\dots,k-1$, and the 
only vectors with lower function value in the $k$-flip neighborhood of a simple vector are other simple vectors. The following is now an immediate corollary of \mbox{Theorem \ref{thm:1}}:
\begin{thm}\label{thm:2}
	Let $k$ be a fixed positive integer. Every $k$-flip local search method applied to $g^{(k)}_n$ starting from $\vx=(0,\dots,0)$ takes exponentially many steps in the size of $g^{(k)}_n$ to reach the local minimum $\vx=(1,\dots,1)$.
\end{thm}

\bibliographystyle{siam}
\bibliography{localsearch_exp_siopt}

\end{document}